\documentclass[12pt,draftcls,onecolumn]{IEEEtran}
\usepackage{epsfig,amssymb,latexsym,color,amsmath,pifont,colordvi,multicol}
\usepackage{subfigure}
\definecolor{backgrey}{rgb}{0.86,0.86,0.86}
\definecolor{dblue}{rgb}{0,0.0,0.5}
\definecolor{dred}{rgb}{0.4,0.2,0}
\definecolor{dgreen}{rgb}{0.0,0.5,0}

\newcommand{\captionfonts}{\small}
\makeatletter  
\long\def\@makecaption#1#2{%
  \vskip\abovecaptionskip
  \sbox\@tempboxa{{\captionfonts #1: #2}}%
  \ifdim \wd\@tempboxa >\hsize
    {\captionfonts #1: #2\par}
  \else
    \hbox to\hsize{\hfil\box\@tempboxa\hfil}%
  \fi
  \vskip\belowcaptionskip}
\makeatother   
\newtheorem{theorem}{Theorem}

\newtheorem{assumption}[theorem]{Assumption}

\newtheorem{remark}[theorem]{Remark}

\newtheorem{corollary}[theorem]{Corollary}

\newtheorem{definition}[theorem]{Definition}

\newtheorem{lemma}[theorem]{Lemma}

\newcommand{\R}{\mathbb{R}}

\title{Synchronization in large-scale nonlinear network systems with uncertain links}

\author{\quad Amit Diwadkar \quad Umesh Vaidya 
\thanks{A. Diwadkar is a Ph.D. student with the Department of Electrical and
Computer Engineering, Iowa State University, Ames, IA, 50011
diwadkar@iastate.edu}
\thanks{U. Vaidya is with the Department of Electrical and
Computer Engineering,
Iowa State University,
Ames, IA, 50011 ugvaidya@iastate.edu}%
}
\begin{document}
\maketitle \thispagestyle{empty} \pagestyle{empty}

\begin{abstract}
In this paper, we study the problem of synchronization with stochastic interaction among network components. The network components dynamics is nonlinear and modeled in Lure form with linear stochastic interaction among network components. To study this problem we first prove the stochastic version of Positive Real Lemma (PRL). The stochastic PRL result is then used to provide sufficient condition for synchronization of stochastic network system. The sufficiency condition for synchronization, is a function of nominal (mean coupling) Laplacian eigenvalues and the statistics of link uncertainty in the form of coefficient of dispersion (CoD). Contrary to the existing literature on network synchronization, our results indicate that both the largest and the second smallest eigenvalue of the nominal Laplacian play an important role in synchronization of stochastic networks. Robust control-based small-gain interpretation is provided for the derived sufficiency condition which allow us to define the margin of synchronization. The margin of synchronization is used to understand the important tradeoff between the component dynamics, network topology, and uncertainty characteristics. For a special class of network system connected over torus topology we provide an analytical expression for the tradeoff between the number of neighbors and the dimension of the torus. Similarly, by exploiting the identical nature of component dynamics computationally efficient sufficient condition independent of network size is provided for general class of network system. Simulation results for network of coupled oscillators with stochastic link uncertainty are presented to verify the developed theoretical framework.
\end{abstract}

\section{Introduction}

The study of network control systems is a topic that has received lots of attention among the research community lately. There is extensive literature on this topic involving both  deterministic and stochastic network systems. Among various problems, the problem of characterizing the stability of estimator and controller design for linear time invariant (LTI) network systems in the presence of channel uncertainty is studied in \cite{ networksystems_foundation_sastry, scl04}. A similar problem involving nonlinear and linear time varying dynamics is studied in \cite{amit_erasure_observation_journal,amit_ltv_journal,Vaidya_erasure_SCL,Vaidya_erasure_stablization}. The results in these papers discover fundamental limitations that arise in the design of stabilizing  controller and estimator in the presence of channel uncertainty.

Another important problem in the study of network systems is that of synchronization of the individual systems interacting over a network with stochastic interactions among network components. Passivity-based tools are used to study the stability problem for deterministic network systems in \cite{Sepul_passivity, pass_network}. Synchronization of interconnected systems from input-output approach has been studied in \cite{arcak_sontag} and shown to have applications in biological networks. These tools provide a systematic procedure for the analysis and synthesis of deterministic network systems. Synchronization of identical nonlinear systems over networks with stochastic link failures was previously studied by the authors in \cite{diwadkar_vaidya_sync}. Master stability function was used to obtain stochastic variational stability in \cite{porfiri_master}. Without assumptions on nonlinearity, the authors were able to provide a necessary condition based on individual system characteristics like Lyapunov exponents and variance of link uncertainty. In this paper, under passivity assumptions on the system dynamics and nonlinearity, we aim to provide a sufficiency condition for synchronization of nonlinear systems over a network with stochastic links. Stability analysis using passivity-based tools for analysis of stochastic systems with additive uncertainty models is studied in \cite{sto_lure1, sto_lure2}. Robustness of synchronized and consensus states to deterministic or stochastic uncertainty has been studied in \cite{trentelman_robust, zampieri_sync, kim_consensus, bassam_coherence, wang_elia}

In this paper, we combine techniques from passivity theory and stochastic systems to provide a sufficient condition for the synchronization of stochastic network systems. Stochastic uncertainty is assumed to enter both multiplicative and additive in system dynamics. We first prove a stochastic version of the Positive Real Lemma (PRL) and provide an Linear Matrix Inequality (LMI)-based verifiable sufficient condition for the mean square exponential stability of stochastic system. An important feature of the stochastic PRL is that the uncertainty enters multiplicatively in the system dynamics. This sufficient condition is then applied to study the problem of synchronization in network of Lur'e systems with  linear but stochastic interactions among the network subsystems. The derived sufficient condition is function of individual component dynamics, network topology captured by the eigenvaues of the nominal Laplacian, and characteristics of stochastic uncertainty captured by coefficient of dispersion (CoD). The CoD  is defined as a ratio of  variance to mean of a random variable and it indicates the amount of clustering behavior in the random variable. A distinguishing feature of the derived sufficient condition is that both the largest and the second smallest (Fiedler eigenvalue) eigenvalues of the nominal Laplacian play a crucial role. The Fielder eigenvalue is well-known in graph theory literature as an indicator of algebraic connectivity of a graph and has been shown to play an important role for problem involving consensus and synchronization in literature. However, to the authors best knowledge, this is the first theoretical demonstration of the role played by the largest eigenvalue of the nominal Laplacian on synchronization, where this eigenvalue plays a role because of the stochastic nature of the network dynamics. Similar influence of the largest eignevalue of the Laplacian is shown computationally in \cite{pecora_carroll_master}. The largest eigenvalue of the nominal Laplacian captures the number of hub nodes and with larger number  of hub nodes in the network uncertainty can propagate faster making it difficult to synchronize. We will elaborate on the significance of the largest eigenvalue of the nominal Laplacian on synchronization in section \ref{section_Lapeigenvalues}. Using the results from robust control theory, small gain theorem-based interpretation is provided for the derived sufficiency condition. The small gain theorem-based interpretation is used to define the margin for synchronization. We provide LMI-based condition for computing the synchronization margin. The synchronization margin play an important role to understand tradeoff between the network component dynamics, network topology, and CoD. In particular for network system with nearest neighbour topology we show that there exists a optimal number of neighbors for which the synchronization margin is largest. This signifies the fact that for network system with uncertainty having too many neighbours is as detrimental to synchronization as having too little neighbours. Furthermore for a special class of torus network \cite{bassam_coherence}, we provide analytical expression to understand the tradeoff between the number of neighbours and the dimension of the torus. Similar results involving limitations and tradeoff for torus networks are derived in \cite{bassam_coherence,wang_elia,xu_ma_torus,xu_elia_torus} for the case of linear time invariant network systems. By exploiting the identical nature of component dynamics we also provide sufficient condition independent of network size for synchronization thereby making the condition attractive from the point of view of computational.

The rest of the paper is structured as follows. In section \ref{prb_def} we formulate the general problem of stabilization of Lur'e systems with parametric uncertainty and prove the main results on the stochastic variant of Positive Real Lemma. The problem of synchronization is formulated and solved using the stochastic variant of PRL in section \ref{sync_prb}. In sections \ref{section_torus_network} and \ref{section_reducedsuffcond}, we discuss synchronization problem on torus networks and provide computationally efficient sufficient condition for synchronization respectively. Simulation results are presented in section \ref{sim_res} followed by conclusions in section \ref{conc}.

\section{Stabilization of Uncertain Lur'e Systems}
In this section, we first present the problem of stochastic stability of a Lur'e system with parametric uncertainty. The uncertainty is modeled as an independent identically distributed (i.i.d.) random processes. The main result of this section proves the stochastic version of the Positive Real Lemma.

\subsection{Problem Formulation}\label{prb_def}
\begin{figure}[ht]
    \centering
    \includegraphics[width=0.3\textwidth]{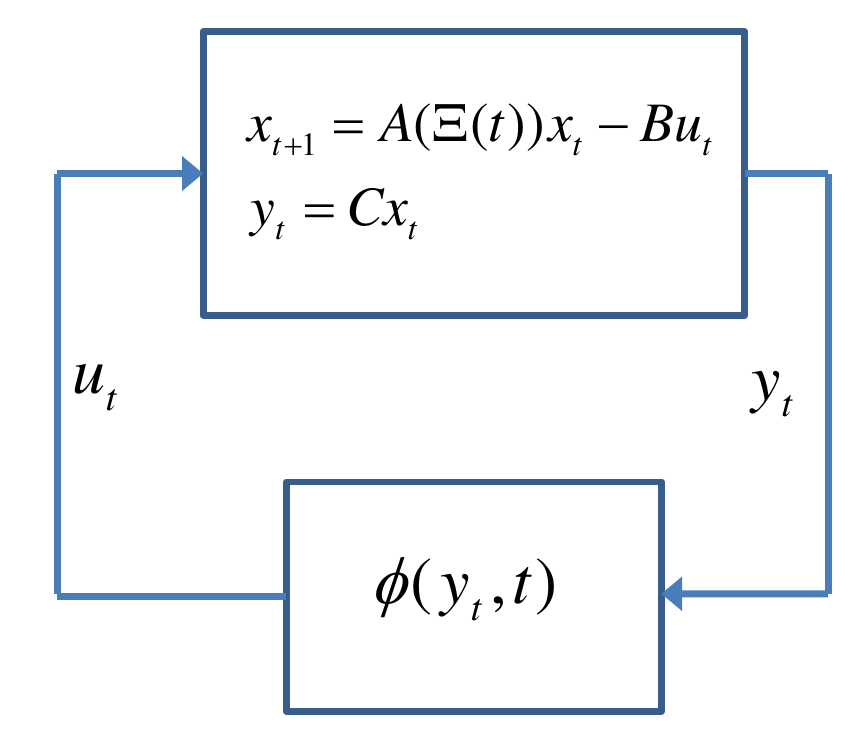}
  \caption{ Schematic of the system with parametric uncertainty.}
  \vspace{-0.25in}
 \label{sche2}
\end{figure}

We consider a Lur'e system, which has parametric uncertainty in the linear system dynamics. The uncertain system dynamics are described as follows:
\begin{align}
x_{t+1} = A \left({\Xi} (t) \right) x_t - B \phi(y_t, t) + v_t, \quad\quad y_t = C x_t \label{unc_lure_sys}
\end{align}
where, $x \in \mathbb{R}^n$, and $y \in \mathbb{R}^m$, $\phi(y_t, t) \in \mathbb{R}^m$ is a nonlinear function, and, $v_t\in \mathbb{R}^n$ is zeros mean additive noise vector with covariance $R_v$. $B \in \mathbb{R}^{n \times m} ~ \text{and} ~ C \in \mathbb{R}^{m \times n}$  are the input and output matrices. The state matrix $A \left({\Xi} (t) \right) \in \mathbb{R}^{n \times n}$ is uncertain, and the uncertainty is characterized by ${\Xi} (t) =[{\xi}_1(t), \dots \xi_M(t)]^T$, where  ${\xi}_i (t)$'s for $i \in \{1,\ldots,M\}$ are i.i.d. random processes with zero mean and variance $\sigma_i^2$, i.e., $E[\xi_i(t)] = 0$, $E[\xi_i(t)^2] = \sigma_i^2$, and $E[\xi_i(t)\xi_j(t)]=0$ for $i\neq j$. The schematic of the system is depicted in Fig. \ref{sche2}.
We make the following assumptions on the nonlinearity $\phi(y_t,t)$
\begin{assumption}\label{nonlin_assume}
The nonlinearity $\phi(y_t,t)$ is a monotonic non-decreasing function of $y_t$ such that, $\phi^\top (y_t,t) \left (y_t- D \phi(y_t,t) \right ) > 0$.
\end{assumption}
The system, described by \eqref{unc_lure_sys}, encompasses a broad class of problems like stabilization under parametric uncertainty, control and observation of Lur'e system over uncertain channel \cite{lure_stab_unc}, and network synchronization of Lur'e systems over uncertain links. Next, we state and prove a stochastic version of the Positive Real Lemma and successively use the result for network synchronization. The stochastic notion of stability that we use is the mean square exponential stability \cite{Hasminskii_book} and is defined as follows:
\begin{definition} \label{def_m_stable}
The system in Eq. \eqref{unc_lure_sys} is mean square exponentially stable if $\exists\; K>0 , ~ \text{and} ~ 0< \beta < 1$, and $L > 0$, such that
\begin{align}
E_{\Xi}  \parallel x_{t} \parallel^2 \le   K  {\beta}^t \parallel x_{0}  \parallel^2 + LR_v, \quad \forall \;x_0  \in \R^{n}. \label{m_st}
\end{align}
 where,  $x_t$ evolves according to \eqref{unc_lure_sys}.
\end{definition}
\begin{remark}
The above definition of mean square exponential stability holds for systems with additive noise. In case the additive noise is absent, the above definitions will reduce to the more familiar definition of mean square exponential stability \cite{Hasminskii_book, lure_stab_unc, amit_erasure_observation_journal}, where $L=0$.
\end{remark}
\subsection{Main Results} \label{prelim}
 The following theorem  is the stochastic version of the Positive Real Lemma providing sufficient condition for the  mean square exponential stability of the stochastic Lur'e system, described by \eqref{unc_lure_sys}.
\begin{theorem} \label{main_kyp1}
Let $\Sigma = D + D^\top $ and $A_T(\Xi(t)) = A(\Xi(t)) - B\Sigma^{-1}C$. Then the uncertain Lur'e system in \eqref{unc_lure_sys} is mean square exponentially stable if -
\begin{enumerate}
\item there exist symmetric positive definite matrices $P$ and $R_P$ such that $\Sigma - B^\top PB > 0$ and,
\begin{align} \label{unc_lure_p}
P =& E_{\Xi(t)} \left [ A_T^\top(\Xi(t)) PA_T(\Xi(t)) + A_T^\top(\Xi(t))PB(\Sigma-B^\top PB)^{-1} B^\top P A_T(\Xi(t))\right ]\nonumber\\
&\quad + C^\top \Sigma^{-1}C+ R_P
\end{align}
\item there exist symmetric positive definite matrices $Q$ and $R_Q$ such that $\Sigma - CQC^\top > 0$ and,
\begin{align} \label{unc_lure_d}
Q =& E_{\Xi(t)} \left [A_T(\Xi(t))QA_T^\top (\Xi(t)) + A_T(\Xi(t))QC^\top (\Sigma-CQC^\top)^{-1} CQA_T^\top (\Xi(t)) \right ]\nonumber  \\
&\quad + R_Q + B^\top \Sigma^{-1}B
\end{align}
\end{enumerate}
\end{theorem}
\begin{IEEEproof} Please refer to the Appendix section for the proof.
\end{IEEEproof}
The generalized version of stochastic Positive Real Lemma, as given by Theorem \ref{main_kyp1}, is now specialized to the case of structured uncertainties. In particular, the structured uncertainties are assumed to be of the form $A(\Xi) = A + \sum_{i=1}^M {\xi}_i A_i $, where $\{\xi_i\}_{i=1}^M$ are zero mean i.i.d. random variables, the mean value having been incorporated in the deterministic part of the matrix given by $A$. The state and output equation for uncertain system becomes,
\begin{align}
x_{t+1} = \left ( A + \sum_{i=1}^M {\xi}_i A_i  \right ) x_t - B \phi(y_t, t) + v_t, \quad\quad y_t = C x_t \label{lin_unc_st}
\end{align}
 The matrices $A_i$, adjoining to the uncertainties, could be pre-determined or could be designed depending on the problem. For instance,
the results developed in \cite{lure_stab_unc} are for the scenario, where the matrix $A_i$ is controller gain. The following Lemma simplifies the generalized stochastic PRL to study the mean square exponential stability of system described by \eqref{lin_unc_st}.
\begin{lemma} \label{new_kyp}
The system, described in \eqref{lin_unc_st}, would be mean square exponentially stable if there exists a symmetric matrix $P > 0$, such that $\Sigma - B^\top PB > 0$ and,
\begin{align} \label{ric_1}
P &= A_0^\top PA_0 + A_0^\top PB(\Sigma-B^\top PB)^{-1}B^\top PA_0 + C^\top \Sigma^{-1}C + R_P  \nonumber \\
&\quad + \sum_{i=1}^M {\sigma}_i^2 \left(A_i^\top PA_i + A_i^\top PB (\Sigma-B^\top PB)^{-1} B^\top PA_i\right)
\end{align}
for some symmetric matrix $R_P>0$ and $A_0 := A - B\Sigma^{-1}C$.
\end{lemma}
\begin{IEEEproof}
We substitute  $\mathcal{A}(\Xi) = A + \sum_{i=1}^M {\xi}_i A_i $ in the \eqref{unc_lure_p} and utilize the fact ${\xi}_i$'s are zero mean i.i.d. random variables with variance ${\sigma}_i^2$. We also $A_T(\Xi) = A + \sum_{i=1}^M{\xi}_i A_i - B\Sigma^{-1}C := A_0 + \sum_{i=1}^M{\xi}_i A_i$. Hence we get,
\begin{align}\label{split_eq1}
& E_{\Xi(t)} \left [A_T^\top (\Xi(t))PA_T(\Xi(t))\right] = A_0^\top PA_0 + \sum_{i=1}^M \sigma_i^2 A_i^\top PA_i
\end{align}
Also we get,
\begin{align}\label{split_eq2}
E_{\Xi(t)} \left [A_T(\Xi(t))^\top PB(\Sigma-B^\top PB)^{-1} B^\top PA(\Xi(t))\right ] &= A_0^\top PB(\Sigma-B^\top PB)^{-1}B^\top PA_0 \nonumber\\
&\quad+ \sum_{i=1}^M \sigma_i^2 A_i^\top PB (\Sigma-B^\top PB)^{-1} B^\top PA_i
\end{align}
Combining equations \eqref{split_eq1} and \eqref{split_eq2} and substituting in \eqref{unc_lure_p} we get the desired result.
\end{IEEEproof}
\begin{corollary} \label{new_kyp_d}
The system, described in \eqref{lin_unc_st}, would be mean square exponentially stable if there exists a symmetric matrix $Q > 0$, such that $\Sigma -CQC^\top  > 0$ and,
\begin{align} \label{ric_2}
Q &= A_0QA_0^\top  + A_0QC^\top(\Sigma-CQC^\top )^{-1}CQA_0^\top  + B^\top\Sigma^{-1}B + R_Q \nonumber \\
& + \sum_{i=1}^M {\sigma}_i^2 \left(A_iQA_i^\top + A_iQC^\top (\Sigma-CQC^\top)^{-1} CQA_i^\top\right)
\end{align}
for some symmetric matrix $R>0$ and $A_0 := A - B\Sigma^{-1}C$.
\end{corollary}
\begin{IEEEproof}
Corollary \ref{new_kyp_d} follows from Theorem \ref{main_kyp1}, Lemma \ref{new_kyp} and duality.
\end{IEEEproof}

%
%
%
%
\section{Synchronization of Lur'e Systems with Uncertain Links}\label{section_main}
In this section, we apply the results developed in the previous section, in analyzing the problem of synchronization of Lur'e systems, coupled through uncertain links. We consider a set of linearly coupled systems in Lur'e form, where the interconnections between these systems, are uncertain in nature. In the subsequent section we derive a sufficiency condition for synchronization over a network, expressed in terms of uncertainty statistics and properties of the mean network, in particular the second smallest and largest eigenvalue of the nominal interconnection Laplacian. The condition could be used to judge whether the coupled system with uncertainty could retain its synchronizability if the links binding the individual subsystems start to fail. Synchronization is achieved if the uncertainty variance satisfies prescribed bounds.
\subsection{Formulation of Synchronization Problem} \label{sync_prb}
 \begin{figure}[t!]
    \centering
    \vspace{-0.25in}
    \includegraphics[width=0.3\textwidth]{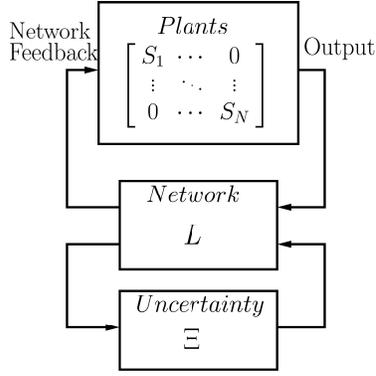}
    \vspace{-0.4in}
  \caption{ Schematic of the interconnected system with uncertain links.}
  \vspace{-0.2in}
 \label{sche1}
\end{figure}
We consider a network of inter-connected systems in Lur'e form.  The individual subsystems are described as follows:,
\begin{equation} \label{component_dynamics}
S_k:=\left\{\begin{array}{ccl}
x^k _{t+1} &=& A  x^k_t -  B {\phi}  (y^k_t, t )\\
y^k_t&=& C  x^k_t,~ ~k=1,\ldots, N
\end{array}\right.
\end{equation}
where, $x^k \in \mathbb{R}^n$, and $y^k \in \mathbb{R}^m$ are the states and the output of $k^{th}$ subsystem. The $\phi(y_n, n) \in \mathbb{R}^l $ is a nonlinear function. The state matrix $A \in  \mathbb{R}^{n \times n}$ is the state matrix for $k^{th}$ subsystem. $B \in \mathbb{R}^{n \times m} ~ \text{and} ~ C\in \mathbb{R}^{m \times n}$  are the input and output matrices of the $k^{th}$ subsystem.  The inter-connected systems interacting with uncertainty through a network are depicted in Fig. \ref{sche1}. The non-linearity satisfies the following assumption,
\begin{assumption} \label{assum_nonlinearity_syn}
The nonlinearity ${\phi}_k(y^k_t,t) \in \R$ is globally Lipschitz monotonically nondecreasing function and $C^1$ function of $y^s_n \in \mathbb{R}$ that satisfies Assumption \ref{nonlin_assume}. Furthermore, it also satisfies the following condition,
\begin{align*}
&\left (  {\phi} \left (y^k_t \right) - {\phi}  \left (y^j_t\right) \right )^\top \left ( (y^k_t - y^j_t)  - D_1 \left ( {\phi} (y^k_t) - {\phi} (y^j_t \right )\right )  > 0,
\end{align*}
for any two systems $S_k$ and $S_j$ and some $\Sigma_1 = D_1 + D_1^\top > 0$.
\end{assumption}
The aforementioned assumption is essential for the synchronization of the network.  Next, we consider coupled subsystems described by equation \eqref{component_dynamics}, that are linearly coupled, and analyze their synchronizability. The coupled system satisfies the following equation,
\begin{align}
x^k_{t+1} = A   x^k_t -  B  {\phi}  \left (y^k_t\right ) + \sum_{j=1}^{N-1} \mu_{kj} G (y^j_t - y^k_t) + v_t, \quad\quad y^k_t = C  x^k_t,~ ~k=1,\ldots, N \label{lap_coup}
\end{align}
where,  $\mu_{kj} \in \mathbb{R}$ represent the coupling link between subsystems $S_k ~\text{and} ~ S_j$, $\mu_{kk} = 0$ and $G \in \mathbb{R}^{n \times m}$.
\begin{remark}
The coupled system as described by \eqref{lap_coup} is the most general form of interaction possible between subsystems. The coupling between subsystems could be either in form of output feedback or state feedback. As the output and states of individual subsystems are related linearly so the form of coupling, as described by \eqref{lap_coup} includes both the output feedback and state feedback.
\end{remark}
Next, we define the graph Laplacian $L_g := \left[ l_{ij}\right] \in \R^{N \times N}$ as following,
\begin{align}
l_{ij} := \mu_{ij}, ~ i \neq j, \;\quad l_{ii} :=-\sum_{j, i \neq j} \mu_{ij}, ~ i=1, \dots N.  \label{laplcn}
\end{align}
Next, all the states of the subsystems are combined to create the states of the coupled system. Finally the coupled system can be rewritten as,
\begin{align}
\tilde{x}_{t+1} = \tilde A \tilde{x}_t - \tilde B \tilde{\phi} \left (\tilde y_t\right ) - \left (L_g \otimes GC \right) \tilde x_t + v_t,\quad\quad \tilde y_t =\tilde C \tilde x_t , \label{coupled_dynamics}
\end{align}
where, $\otimes $ is the Kronecker product, $I_n$ is an $n\times n$ Identity matrix and,
\begin{equation*}
\tilde A:=I_N\otimes A=
\begin{bmatrix}
A & 0 & \cdots & 0 \\
0 & A & \cdots & 0 \\
\vdots & \vdots & \ddots & \vdots \\
0 & 0 & \cdots & A
\end{bmatrix}
\end{equation*}
We similarly define $\tilde{B} := I_N\otimes B$, $\tilde{C} := I_N\otimes C$, $\tilde{D}_1 := I_N\otimes D_1$ and $\tilde{\Sigma}_1 := \tilde D_1 + \tilde D_1^\top > 0$. We also define
$\tilde x_t  = [ (x^1_t)^\top \ldots (x^N_t)^\top  ]^\top, \;
\tilde y_t  = [ (y^1_t)^\top \ldots (y^N_t)^\top  ]^\top, \;
\tilde {\phi}_t = [ ({\phi}^1_t)^\top\ldots ({\phi}^N_t)^\top  ]^\top$.

%
\subsection{Modeling Uncertain Links}\label{main}
We are now ready to study the problem of synchronization where the links of the graph are uncertain ( i.e. entries of the Laplacian matrix are uncertain). Let
\[ E_U=\{(i,j)| \text{the link}\; (i, j) \;\text{is uncertain}, ~ i > j \} \]
be the collection of uncertain links in the network. Hence, for links $(i,j) \in E_U$, we have $l_{ij} = \mu_{ij} + \xi_{ij}, \; i\neq j$, where $\xi_{ij}$ are zero mean i.i.d. random variables with variance $\sigma_{ij}^2$. If $(i,j) \notin E_U$ when we have $l_{ij} = \mu_{ij}, \; i\neq j$ to be purely deterministic as in the previous subsection. This framework allows us to study synchronization for Lur'e type systems with a deterministic weighted Laplacian as a special case. Let $\Xi = \{\xi_{ij}\}_{(i,j)\in E_U}$. Then, the uncertain graph Laplacian $L_g(\Xi)$ will be given as,
\begin{align}
L_g(\Xi) = L + \sum_{(i,j)\in E_U}\xi_{ij}L_{ij}
\end{align}
where $L$ is the nominal part of the uncertain graph Laplacian $L_g(\Xi)$, which may be written as $L = L_d + L_u$, where $L_d$ is the part of the Laplacian constructed from $\mu_{ij}$ for purely deterministic edges $(i,j) \notin E_U$, while $L_u$ is constructed from the mean weights $\mu_{ij}$ for the uncertain edges $(i,j) \in E_U$. We may also write $L_{ij} = \ell_{ij}\ell_{ij}^\top$ where $\ell_{ij} := [\ell_{ij}(1),\ldots,\ell_{ij}(N)]^\top \in \mathbb{R}^N$ is a column vector given by
\begin{align*}
\ell_{ij}(k) =
\left\{\begin{array}{cc}
0 &\text{if}\; k \neq i \neq j \\
 1 &\text{if}\; k=i\\
 -1 &\text{if}\; k=j
\end{array}\right.
\end{align*}
We are interested in finding a sufficiency condition involving $\sigma_{ij}^2$ for $(i,j)\in E_U$, which would guarantee the mean square exponential synchronization. The coupled network of Lur'e system can be written as,
\begin{align}
\tilde{x}_{t+1} = \left( \tilde A - \left( L_g(\Xi) \otimes G C \right)\right ) \tilde{x}_n - \tilde B \tilde{\phi} \left (\tilde y_t\right ) + v_t,\quad\quad \tilde y_t =\tilde C \tilde x_t \label{coup_dyn_comp}
\end{align}
We would analyze the stochastic synchronization of system, described by \eqref {coup_dyn_comp}. We start with following definition of mean square exponential synchronization.
\begin{definition} \label{sto_syn_def}
The system, described by \eqref{coup_dyn_comp} is mean square exponentially synchronizing if there exists a $\beta < 1$, $K(\tilde e_0) > 0 $, and, $L > 0$  such that,
\begin{eqnarray}
& E_{\Xi} \parallel x^k_{t} - x^j_{t}\parallel^2 \leq   \bar K(\tilde e_0) {\beta}^t \parallel x^k_{0} - x^j_{0}\ \parallel^2 + LR_v, \;\;\;\forall k,j \in [1,N] \label{mss_sync}
\end{eqnarray}
where, $\tilde e_0$ is function of difference $\parallel x_t^i-x_t^\ell\parallel^2$ for $i,\ell \in [1,N]$ and $\bar K(0)=K$ for some constant $K$.
\end{definition}

We now apply change of coordinates to decompose the system dynamics on and off the synchronization subspace. The synchronization subspace is given by ${\bf 1}=[1,\ldots,1]^\top$. We show that the dynamics on the synchronization subspace is decoupled from the dynamics off the manifold and is essentially described by the dynamics of the individual system. The dynamics on the synchronization subspace itself could be stable, oscillatory, or complex.
Let $L_m = V_m\Lambda_mV_m^\top$ where $V_m$ is an orthonormal set of vectors given by $V_m = \left[\frac{\bf{1}}{\sqrt{N}}\; U_m\right]$, ${\bf 1} = [1 \;\cdots\;1]^\top$ and $U_m$ is orthonormal set of vectors also orthonormal to ${\bf 1}$. Let $\tilde z_t = \left(V_m^\top\otimes I_n\right)\tilde x_t$. Multiplying (\ref{coup_dyn_comp}) from the left by $V_m^\top\otimes I_n$ we get
\begin{align}
\tilde{z}_{t+1} &= \left( \tilde A - \left( V_m^\top L_g(\Xi)V_m \otimes G C \right)\right ) \tilde{z}_t - \tilde B \tilde{\psi} \left (\tilde w_t\right ) + \vartheta_t \label{coup_dyn_comp1}
\end{align}
where $\tilde w_t =\tilde C \tilde z_t$, $\tilde{\psi}_t = \left(V_m^\top \otimes I_n\right)\tilde{\phi}\left(\tilde y_t\right)$, and $\upsilon_t = \left(V_m^\top \otimes I_n\right)v_t$. We can now write
\begin{align}
\tilde z_t = \left[\begin{array}{cc}\bar{x}_t^\top & \hat{z}_t^\top \end{array}\right]^\top,\quad
\tilde{\psi}_t= \left[\begin{array}{cc}
\bar{\phi}_t^\top & \hat{\psi}_t^\top
\end{array}\right]^\top\label{split_vec}
\end{align}
where
\begin{align}\label{transform1}
\bar{x}_t &:= \frac{\bf{1}}{\sqrt{N}}\tilde{x}_t = \frac{1}{\sqrt{N}}\sum_{k = 1}^N x^k_t,\quad \hat{z}_t := \left(U_m^\top \otimes I_n\right)\tilde{x}_t \\
\bar{\phi}_t &:= \frac{\bf{1}}{\sqrt{N}}\tilde{\phi}\left(\tilde y_t\right) = \frac{1}{\sqrt{N}}\sum_{k = 1}^N \phi(y^k_t),\quad \hat{\psi}_t := \left(U_m^\top \otimes I_n\right)\tilde{\phi}\left(\tilde y_t\right)
\end{align}
Substituting \eqref{split_vec} in \eqref{coup_dyn_comp1} we get
\begin{align}
\bar{x}_{t+1} &=  A \bar{x}_t - B \bar{\phi} \left (\bar y_t\right ) + \bar{v}_t \nonumber\\
\hat{z}_{t+1} &= \left( \hat A -  \left( U_m^\top L_g(\Xi) U_m \otimes G C \right)\right ) \hat{z}_t - \hat B \hat{\psi} \left (\hat w_t\right ) +\hat{\vartheta}_t \label{coup_dyn_comp2}
\end{align}
where $\hat w_t =\hat C \hat z_t$,
$\hat A:=I_{N-1}\otimes A,\; \hat B:=I_{N-1}\otimes B,\;
\hat C:=I_{N-1}\otimes C$, and $\hat D_1:=I_{N-1}\otimes D_1$.
We now show that for the synchronization of system (\ref{coup_dyn_comp}), we only need to stabilize $\hat{z}_t$ dynamics. The stability of the system with state $\hat z_t$, implies the synchronization of the actual coupled system. This feature is exploited to derive sufficiency condition for stochastic synchronization of the coupled system. In the following Lemma we show the connection between the stability of the described by \eqref{coup_dyn_comp2} to the synchronization of the system described by \eqref {coup_dyn_comp}.

\begin{lemma} \label{con_syn_sta}
Mean square exponential stability of system described by \eqref {coup_dyn_comp2} implies mean square exponential synchronization of the system \eqref{coup_dyn_comp} as given by Definition \ref{sto_syn_def}.
\end{lemma}
Please refer to the Appendix section of this paper for the proof.
%
In the following subsection we will provide sufficiency conditions for the mean square exponential synchronization of \eqref{coup_dyn_comp} by proving sufficiency conditions for mean square exponential stability of \eqref{coup_dyn_comp2}. But first, we rewrite the equation \eqref{coup_dyn_comp2} in a more suitable format. We note that $L_g(\Xi) =~ L_m + \sum_{E_U}\xi_{ij}L_{ij}$, and $L_m = V_m\Lambda_mV_m^\top $ where $V_m = \left[\frac{\bf{1}}{\sqrt{N}}\; U_m\right]$. Hence we have
\begin{align*}
U_m^\top L_g(\Xi)U_m = U_m^\top L_mU_m + \sum_{E_U}\xi_{ij}U_m^\top L_{ij}U_m
:= \hat{\Lambda}_m + \sum_{E_U}\xi_{ij}\hat{\ell}_{ij}\hat{\ell}_{ij}^\top
\end{align*}
where $L_{ij} = \ell_{ij}\ell_{ij}^\top $, $\hat{\ell}_{ij} = U_m^\top \ell_{ij}$ and $\hat{\Lambda}_m := U_m^\top L_mU_m$ such that
\begin{align*}
\Lambda_m = V_m^\top L_mV_m = \left[\begin{array}{cc}
0 & 0\\
0 & U_m^\top L_mU_m
\end{array}\right] = \left[\begin{array}{cc}
0 & 0\\
0 & \hat{\Lambda}_m
\end{array}\right]
\end{align*}
Let $\mathcal{I} = \{\alpha_k\}_{k=1}^M$, $M = |E_U|$ be an indexing on uncertain edges in $E_U$. If index $\alpha_k$ corresponds to edge $(i,j) \in E_U$ then let $A_{\alpha_k} := U_m^\top L_{ij}U_m\otimes GC = \hat{\ell}_{ij}\hat{\ell}_{ij}^\top  \otimes GC$. Thus we can write equation \eqref{coup_dyn_comp2} as
\begin{align}
\hat z_{t+1} = \left(\hat{A} - \hat{\Lambda}_m\otimes GC - \sum_{\alpha_k\in\mathcal{I}}\xi_{\alpha_k}A_{\alpha_k}\right)\hat{z}_t - \hat{B}\hat{\psi}_t(\hat{w}_t) + \hat{\vartheta}_t
\label{sync_dyn}
\end{align}
\subsection{Sufficiency Condition for Synchronization with Uncertain Links}
In previous subsection, we have shown that mean square exponential stability of \eqref{sync_dyn} guarantees the mean square exponential synchronization of the coupled network of Lur'e system as given by \eqref{coup_dyn_comp}. In the preceeding section, we have derived sufficiency condition for mean square exponential stability of Lur'e system. In this subsection, we combine these two results to obtain sufficiency condition for mean square exponential synchronization of the network of Lur'e systems. The following Lemma provides the sufficiency condition for mean square exponential synchronization.

\begin{theorem} \label{suf_syn1}
The system described by \eqref{coup_dyn_comp} is mean square exponential synchronizing if there exists a symmetric positive definite matrix $\mathcal{P} \in \mathbb{R}^{(N-1)n \times (N-1)n}$ such that,
\begin{align} \label{ric_syn}
\mathcal{P} =& (\hat{A}_0 - \hat\Lambda_m\otimes GC)^\top \mathcal{P}(\hat{A} - \hat\Lambda_m\otimes GC) + \sum_{\mathcal{I}} \sigma_{\alpha_k}^2A_{\alpha_k} ^\top  \mathcal{P} A_{\alpha_k}\nonumber \\
& + (\hat{A}_0-\hat\Lambda_m\otimes GC)^\top \mathcal{P}\hat{B}\left(\hat \Sigma_1 - \hat B^\top  \mathcal{P} \hat B\right)^{-1}\hat B^\top \mathcal{P}(\hat{A}_0-\hat\Lambda_m\otimes GC) \nonumber \\
&+ \sum_{\mathcal{I}} {\sigma}_{\alpha_k}^2 A_{\alpha_k}^\top  \mathcal{P}\bar B \left(\hat \Sigma_1-\hat B^\top  \mathcal{P} \hat B\right)^{-1} \hat B^\top \mathcal{P} A_{\alpha_k}  + \mathcal{R}
\end{align}
and $\hat \Sigma_1-\hat B^\top \mathcal{P} \hat B > 0$ for some symmetric matrix $\mathcal{R} > 0$ and $\hat{A}_0 := \hat{A} - \hat{B}\hat\Sigma_1^{-1}\hat{C} = I_{N-1}\otimes A_0$, $A_0 = A - B\Sigma_1^{-1}C$.
\end{theorem}

\begin{IEEEproof}
The proof follows from \eqref{coup_dyn_comp}, \eqref{sync_dyn}, Lemma \ref{con_syn_sta} and Theorem \ref{main_kyp1}.
\end{IEEEproof}

\subsection{Small Gain Theorem Based Interpretation}

In this subsection, we provide a Small Gain Theorem based interpretation of the sufficient condition for mean square exponential synchronization, as provided in Theorem \ref{suf_syn1}. In Theorem \ref{suf_syn1}, we have derived the sufficient condition for mean square exponential synchronization of coupled $n$-dimensional Lur'e systems with multiple link uncertainties. To provide the Small Gain based interpretation for the derived sufficiency condition we will make following assumption on the uncertainty statistics
\begin{assumption}\label{assume_smallgain} We assume that the all the stochastic interaction uncertainties have identical variance i.e., $E[\xi_i^2(t)]=\sigma^2$ for $i=1,\ldots,m$. Furthermore, we assume that the additive noise term $v_t\equiv 0$.
\end{assumption}
Using assumption \ref{assume_smallgain}, we can rewrite the sufficiency condition as follows:
\begin{align} \label{suff_ric1}
\mathcal{P} > &  (\hat{A}_0 - \hat\Lambda_m\otimes GC)^\top \mathcal{P}(\hat{A} - \hat\Lambda_m\otimes GC) + \sigma^2\sum_{\mathcal{I}} A_{\alpha_k} ^\top  \mathcal{P} A_{\alpha_k}\nonumber \\
& + (\hat{A}_0-\hat\Lambda_m\otimes GC)^\top \mathcal{P}\hat{B}\left(\hat \Sigma_1 - \hat B^\top  \mathcal{P} \hat B\right)^{-1}\hat B^\top \mathcal{P}(\hat{A}_0-\hat\Lambda_m\otimes GC) \nonumber \\
& +  \sigma^2\sum_{\mathcal{I}} A_{\alpha_k}^\top  \mathcal{P}\bar B \left(\hat \Sigma_1-\hat B^\top \mathcal{P} \hat B\right)^{-1} \hat B^\top \mathcal{P} A_{\alpha_k}
\end{align}

The system given by Eq. \eqref{suff_ric1}, can further be written in the following input-ouput form,
\begin{align}\label{sys_inout}
\hat{z}_{t+1} &= \left(\hat{A}_0-\hat{\Lambda_m}\otimes GC\right)\hat{z}_t - \hat{\Upsilon}\otimes G \hat{\eta}_t + \hat{B}\hat{\nu}_t, \\
\hat{w}_t &= \hat{C}\hat{z}_t, \quad \hat{\omega}_t = \hat{\Upsilon}^\top \otimes C\hat{z}_t,\\
\hat{\nu}_t &= \hat\Sigma_1 \hat{w}_t - \hat{\psi}(\hat{w}_t), \quad \hat{\eta}_t = \Xi \hat{\omega}_t,
\end{align}
where $\Xi = \text{diag}\{\xi_1,\ldots,\xi_M\}$, $E[\xi_k] = \sigma^2$, and $\hat{\Upsilon} = [\hat{\ell}_1\quad \hat{\ell}_2\quad \cdots \hat{\ell}_M]$.

\begin{figure}[ht]
\begin{center}
\includegraphics[width=4.5in]{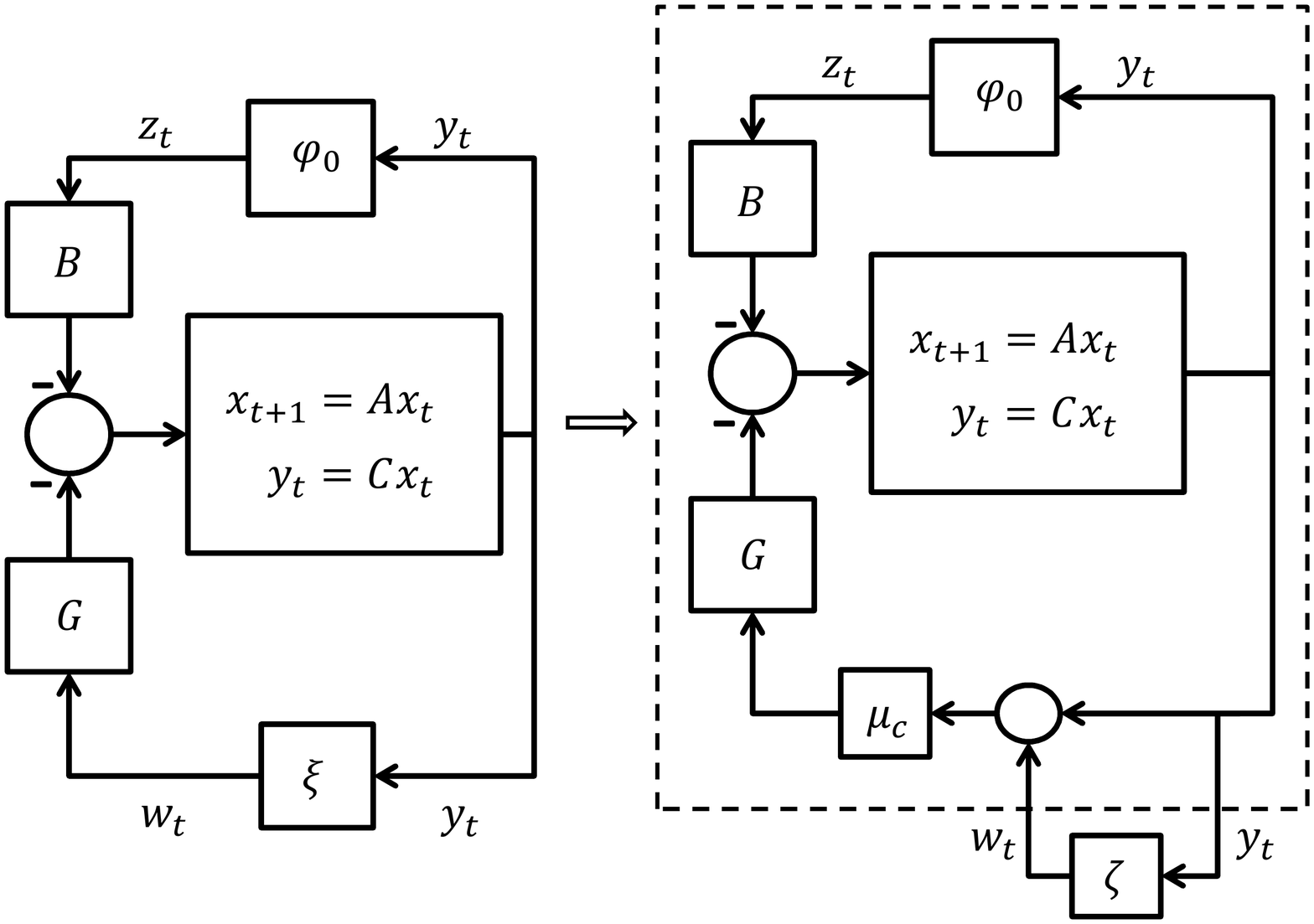}
  \caption{Schematic digram of system with stochastic uncertainty feedback loop and nonlinearity feedback loop}
\vspace{-0.2in}
 \label{diag}
 \end{center}
\end{figure}

This can be represented in a schematic diagram as shown in Fig. \ref{diag}. We will now try to interpret the mean square exponential synchronization condition in Eq. \eqref{suff_ric1}, as a loop gain stability condition for the mean deterministic part and the stochastic uncertainty from output $y_t$ to input $w_t$, as shown in Fig. \ref{diag}. This is similar to the Small Gain interpretation in robust control theory \cite{paganini_dull} or stochastic robust control theory \cite{scl04}. We now present a theorem which illustrates the Small Gain like nature of the sufficient condition.

\begin{theorem}\label{sm_gain_thm}
The input-output system described in Eq. \eqref{sys_inout} is mean square exponentially stable if
\begin{align}\label{sm_gn}
1 > \sigma^2 \rho(M)^2,
\end{align}
where,
\begin{align}\label{M_mat}
M = \begin{bmatrix} \hat{S}^{-\frac{1}{2}}A_{\alpha_1}\left(\hat{T} - \hat A_0^\top \hat S^{-1}\hat A_0\right)^{-\frac{1}{2}} \\ \vdots \\ \hat{S}^{-\frac{1}{2}}A_{\alpha_M}\left(\hat{T} - \hat A_0^\top\hat S^{-1}\hat A_0\right)^{-\frac{1}{2}} \end{bmatrix},
\end{align}
$\hat S=\mathcal P^{-1} - \hat B\hat \Sigma_1^{-1}\hat B^\top$, and, $\hat T=\mathcal P - \hat C^\top\hat \Sigma_1^{-1}\hat C$.
\end{theorem}
\begin{IEEEproof}
Defining $\hat S=\mathcal P^{-1} - \hat B\hat \Sigma_1^{-1}\hat B^\top$, and, $\hat T=\mathcal P - \hat C^\top \hat \Sigma_1^{-1}\hat C$, and applying Matrix Inversion Lemma to Eq. \eqref{suff_ric1}, we obtain,
\begin{align}\label{compact_suff}
\hat T > \hat A_0^\top \hat S^{-1}\hat A_0 + \sigma^2\sum_{\mathcal{I}}A_{\alpha_k}^\top \hat S^{-1}A_{\alpha_k}.
\end{align}
Since $\hat S$ and $\hat T - \hat A_0^\top \hat S^{-1}\hat A_0$ are positive definite matrices, the can be written as a square of positive definite matrices $\hat S^{-\frac{1}{2}}$ and $(\hat T - \hat A_0^\top \hat S^{-1}\hat A_0)^{-\frac{1}{2}}$, respectively. Hence, from \eqref{compact_suff} we obtain,
\begin{align}\label{compact1}
I > \sigma^2 \sum_{\mathcal{I}}\left(\hat T - \hat A_0^\top \hat S^{-1}\hat A_0\right)^{-\frac{1}{2}}A_{\alpha_k}^\top \hat S^{-1}A_{\alpha_k}\left(\hat T - \hat A_0^\top \hat S^{-1}\hat A_0\right)^{-\frac{1}{2}}.
\end{align}
Defining $M$ as given in \eqref{M_mat} 
 and substituting in \eqref{compact1}, we obtain,
\begin{align}\label{compact2}
I > \sigma^2M^\top M.
\end{align}
A sufficient condition for \eqref{compact2} is given by
\begin{align}
1 > \sigma^2\rho(M)^2,
\end{align}
where, $\rho(M)$ is the singular value of the matrix $M$.
\end{IEEEproof}

\begin{remark}
We observe the condition \eqref{sm_gn} provided in Theorem \ref{sm_gain_thm}, is similar to the small gain condition for stochastic stability, where $\rho(M)^2$, can be considered as the mean square gain of the mean deterministic system. The system norm can be computed using an iterative solving of the Riccati for mean square exponential stability given in \eqref{suff_ric1} by writing it as an LMI.
\end{remark}

The above sufficiency condition is very difficult to verify for large size networks due to computational complexity associated with solving the Riccati equation. In particular the matrix $\cal P$ is of size $(N-1)n\times (N-1)n$ having $\frac{(N-1)^2n^2 + (N-1)n}{2}$ variables to be determined. The number of variables increases quadratically with change in system dimension or size of network. In the following we will discuss two different  scenarios which will allow us to reduce the computational complexity. In particular, in section \ref{section_torus_network} we study the case where the nonlinear component dynamics is connected over torus network. The analytical formula for the eigenvalues of the nominal Laplacian with torus geometry will substantially reduce the computational complexity. In section \ref{section_reducedsuffcond}, we exploit the identical nature of system dynamics to provide more conservative sufficient condition but with substantially reduced computational efforts. The sufficiency condition is based upon a single representative dynamical system modified using network characteristics, thereby making it independent of network size.

\section{Spatially Periodic Nonlinear Networks}\label{section_torus_network}
In this subsection we study synchronization of Lure system connected over spatially periodic networks. Spatially periodic networks or torus networks are studied in the context of linear time invariant system. Issues related to fundamental limitations for coherency in consensus network were addressed in \cite{bassam_coherence}. Similarly the problem of consensus over torus network with stochastic interaction uncertainty with LTI dynamics is studied in \cite{wang_elia, xu_ma_torus}. One of the important characteristics of torus network which helps in simplifying their analysis is that their topological properties like the nominal Laplacian eigenvalues has an analytical expression. A simplest torus network is given by a nearest neighbor network, where each agent is connected to adjacent two neighbors. A higher dimensional torus network can be constructed by adding identical nearest neighbor networks along each dimension. Consider a nearest neighbor network with $d$ dimensions, and $N$ agents with $k$ nearest neighbors in each dimension, having a Laplacian matrix given by $L_{N,k,d}$. If a $d+1$ dimensional $N$-neighbor network is constructed with Laplacian $L_{N,k,d+1}$ then we have,
\begin{align}\label{lap_torus}
L_{N,k,d} = I_N\otimes L_{N,k,d-1} + L_{N,k,1}\otimes I_{N^{d-1}} = \sum_{i=0}^{d-1} I_{N^{d-1-i}}\otimes L_{N,k,1} \otimes I_{N^i}
\end{align}
where $\otimes$ denotes the Kronecker product of matrices. Using the eigenvalue property of Kronecker products of matrices, we obtain
\begin{align}\label{eig_torus}
\lambda(L_{N,k,d}) = \sum_{i = 1}^d \lambda_{k_i}(L_{N,k,1})
\end{align}
where $\lambda_{k_i}(L_{N,k,1})$ are the eigenvalues of the 1-d torus with Laplacian $L_{N,k,1}$, for all $k_i \in \{1,2,\ldots,N\}$. Suppose the smallest non-zero eigenvalue of $L_{N,k,1}$ is $\lambda_2$ and largest eigenvalue is $\lambda_N$, then the smallest non-zero eigenvalue of $L_{N,k,d}$ is given by $\tilde\lambda_2 := \lambda_2(L_{N,k,d}) = \lambda_2$, and the largest eigenvalue of $L_{N,k,d}$ is given by $\tilde\lambda_N := \lambda_N(L_{N,k,d}) = d\lambda_N$. We will use these results to prove the results in this section.

We now consider a system \eqref{coup_dyn_comp} on a $d$-torus with $N$ agents and $k$ neighbors along each dimension of the torus network. We also assume that all the links in the network are uncertain and all the links have identical weight $\mu$ and are affected by the same zero mean uncertainty $\xi$ with variance $\sigma2$. This allows us to write \eqref{coup_dyn_comp} as follows
\begin{align}\label{torus_sys}
\tilde{x}_{t+1} = \left( \tilde A - (\mu + \xi)\left( L_{N,k,d} \otimes G C \right)\right ) \tilde{x}_n - \tilde B \tilde{\phi} \left (\tilde y_t\right ) + v_t,\quad\quad \tilde y_t =\tilde C \tilde x_t
\end{align}
Let the eigenvectors of $L_{N,k,d}$ be given by $V_g$ and the diagonal matrix of eigenvalues be given by $\Lambda_{N,k,d}$. Hence we obtain $V_g'L_{N,k,d}V_g = \Lambda_{N,k,d}$. Applying the transformation $V'\otimes I_n$ to $\tilde{x}_t$, we obtain,
\begin{align}\label{torus_sys1}
\tilde{z}_{t+1} = \left( \tilde A - (\mu + \xi)\left( \Lambda_{N,k,d} \otimes G C \right)\right ) \tilde{z}_t - \tilde B \tilde{\phi} \left (\tilde y_t\right ) + \tilde{v}_t,\quad\quad \tilde y_t =\tilde C \tilde z_t
\end{align}
Rewriting the mean dynamics separately as $\hat{x}_t$, from \eqref{torus_sys1} we obtain,
\begin{align}\label{torus_sys2}
\bar{x}_{t+1} &= A\bar{x}_t - B\bar{\phi}(\bar{y}_t) + \bar{v}_t,\quad \quad \bar{y}_t = C\bar{x}_t\\
\hat{z}_{t+1} &= \left( \hat A - (\mu + \xi)\left( \hat\Lambda_{N,k,d} \otimes G C \right)\right ) \hat{z}_t - \hat B \hat{\phi} \left (\hat y_t\right ) + \hat{v}_t,\quad\quad \hat y_t =\hat C \hat z_t
\end{align}
We invoke Lemma \ref{suf_syn1} to obtain the stability condition for the system $\hat{z}_t$ as given in \eqref{torus_sys2}. Using Lemma \ref{suf_syn1} we can provide the following lemma for the mean square exponential stability of $\tilde{x}_t$.

\begin{lemma}\label{torus_suff_condn}
The system of agents interacting over a $d$-dimensional torus network as given in \eqref{torus_sys}, is mean square exponentially synchronizing if there exist positive definite matrices $P_i > 0$ for $i \in \{2,3,\ldots,N^d\}$, such that $\Sigma_1 -B'P_iB > 0$ for all $i$, and
\begin{align}\label{suff_condn_tori}
P_i &> (A_0-\mu\tilde\lambda_iGC)'P_i(A_0-\mu\tilde\lambda_iGC) + (A_0-\mu\tilde\lambda_iGC)'P_iB\left(\Sigma_1 - B'P_iB\right)^{-1}B'P_i(A_0-\mu\tilde\lambda_iGC) \nonumber\\
& \qquad \sigma^2\lambda_i^2\left(C'G'P_iGC + C'G'P_iB\left(\Sigma_1 - B'P_iB\right)^{-1}B'P_iGC\right) + C'\Sigma_1C,
\end{align}
where $A_0 = A - B\Sigma_1^{-1}C$. Furthermore, the condition \eqref{suff_condn_tori} exists if the follwoing linear matrix inequality (LMI) is satisfied for all $P_i > 0$,
\begin{align}\label{tori_lmi}
\begin{bmatrix}
P_i & (A_0 - \mu\tilde\lambda_iGC)'P_i & \sigma\tilde\lambda_iC'G'P_i & (A_0 - \mu\tilde\lambda_iGC)'P_iB & \sigma\tilde\lambda_iC'G'P_iB \\
P_i(A_0 - \mu\tilde\lambda_iGC) & P_i & 0 & 0 & 0\\
\sigma\tilde\lambda_iP_iGC & 0 & P_i & 0 & 0\\
B'P_i(A_0 - \mu\tilde\lambda_iGC) & 0 & 0 & \Sigma_1 - B'P_iB & 0\\
\sigma\tilde\lambda_iB'P_iGC & 0 & 0 & 0 & \Sigma_1 - B'P_iB
\end{bmatrix} > 0
\end{align}
\end{lemma}

\begin{IEEEproof}The proof follows from Lemma \ref{suf_syn1} and the fact that $\hat{z}$ is a set of uncertain decoupled equations.
\end{IEEEproof}

The above LMI in \eqref{tori_lmi} or Riccati equation in \eqref{suff_condn_tori} is difficult to solve for higher dimensional systems as you have to solve that for and $n\times n$ matrix $P_i$ for all possible eigenvalues. In the following theorem we study the above condition for simple scalar systems with the assumptions $A = a$, $B = 1$, $C = 1$ and $D = \frac{\delta}{2}$ with dynamics similar to \eqref{torus_sys} given by,
\begin{align}\label{torus_sys_1d}
\tilde{x}_{t+1} = \left( aI - (\mu + \xi)g L_{N,k,d}\right ) \tilde{x}_t - \tilde{\phi} \left (\tilde x_t\right ) + v_t.
\end{align}
Hence for scalar agents in teracting over torus networks, the condition for mean square exponential stability is given by the following theorem.

\begin{theorem}\label{1Dtorus_suff_condn}
The system of scalar agents interacting over a $d$-dimensional torus network as given in \eqref{torus_sys_1d}, is mean square exponentially synschronizing if there exist positive definite scalars $\delta > p_i > 0$ for $i \in \{2,3,\ldots,N^d\}$, such that $\delta > p_i$ for all $i$, and
\begin{align}\label{1D_suff_tori}
p_i &> \frac{(a_0-\mu\tilde\lambda_ig)^2\delta p_i}{\delta - p_i} + \tilde\lambda_i^2\sigma^2\frac{g^2\delta p_i}{\delta - p_i} + \frac{1}{\delta},
\end{align}
where $a_0 = a - \frac{1}{\delta}$. The conditions in \eqref{suff_condn_tori} are satisfied if and only if
\begin{align}\label{1D_suff_tori_margin}
\left(1 - \frac{1}{\delta}\right)^2 > \max\{\alpha_2^2,\alpha_{N^d}^2\},
\end{align}
where $\alpha_i^2 := (a_0-\mu\tilde\lambda_ig)^2 + \sigma^2\tilde\lambda_i^2g^2$, for all $i \in \{2,\ldots,N^d\}$. Therefore, we can define the margin of synchronization for a given variance of uncertainty as
\begin{align}\label{sync_margin}
\rho_{SM} := 1 - \sigma^2 \left(\frac{\tilde{\lambda}_{sup}^2g^2}{\left(1 - \frac{1}{\delta}\right)^2 - \left(a_0 - \mu\tilde{\lambda}_{sup}g\right)^2}\right),
\end{align}
where $\tilde{\lambda}_{sup} = \arg\max_{\tilde{\lambda_2},\tilde{\lambda}_{N^d}}\alpha_i^2$.
\end{theorem}

\begin{IEEEproof}
From Lemma \ref{torus_suff_condn} we know there exists scalars $\delta > p_i > 0$ such that
\begin{align}\label{1D_suff_tori2}
p_i &> \frac{(a_0-\mu\tilde\lambda_ig)^2\delta p_i}{\delta - p_i} + \tilde\lambda_i^2\sigma^2\frac{g^2\delta p_i}{\delta - p_i} + \frac{1}{\delta} > \frac{\delta p_i}{\delta - p_i}\left( (a_0-\mu\tilde\lambda_ig)^2 + \tilde{\lambda}_i^2\sigma^2g^2 \right) + \frac{1}{\delta}
\end{align}
We can rewrite Eq. \eqref{1D_suff_tori2} as
\begin{align}\label{poly_suff}
\left(p_i - \frac{1}{\delta} \right)\left( \frac{1}{p_i} - \frac{1}{\delta} \right) > \alpha_i^2,
\end{align}
where we define $\alpha_i^2 := (a_0-\mu\tilde\lambda_ig)^2 + \tilde{\lambda}_i^2\sigma^2g^2$. Now using the AM-GM inequality we can write
\begin{align}\label{am_gm_ineq}
\left(p_i - \frac{1}{\delta} \right)\left( \frac{1}{p_i} - \frac{1}{\delta} \right) &= 1 -\left(p_i + \frac{1}{p_i}\right)\frac{1}{\delta} + \frac{1}{\delta^2},\nonumber\\
&< 1 - \frac{2}{\delta} + \frac{1}{\delta^2},\nonumber\\
&< \left(1 - \frac{1}{\delta}\right)^2.
\end{align}
Now using \eqref{am_gm_ineq} in \eqref{poly_suff}, we obtain a sufficient condition for \eqref{poly_suff} given by
\begin{align}\label{poly_suff_1}
\left(1 - \frac{1}{\delta}\right)^2 > \alpha_i^2,
\end{align}
for all $i \in \{2,\ldots,N^d\}$. Now suppose \eqref{poly_suff_1} is true then there exists $\epsilon_i > 0$ such that,
\begin{align}\label{eps_diff}
\left(1 - \frac{1}{\delta}\right)^2 - \alpha_i^2 = 2\frac{\epsilon_i^2}{1 + \epsilon_i}.
\end{align}
Taking $p_i = 1 + \epsilon_i$ we obtain
\begin{align}\label{p_delta}
\left( 1 - \frac{1}{\delta}\right)^2 - \frac{\epsilon_i^2}{1 + \epsilon_i} &= 1 - \frac{2}{\delta} + \frac{1}{\delta^2} - \frac{\epsilon_i^2}{1 + \epsilon_i},\nonumber\\
&= 1 - \left(1 + \epsilon_i +\frac{1}{1+\epsilon_i}\right)\frac{1}{\delta} + \frac{1}{\delta^2},\nonumber\\
&= 1 - \left(p_i +\frac{1}{p_i}\right)\frac{1}{\delta} + \frac{1}{\delta^2},\nonumber\\
&= \left(p_i - \frac{1}{\delta}\right)\left(\frac{1}{p_i} - \frac{1}{\delta}\right).
\end{align}
Substituting \eqref{p_delta} in \eqref{eps_diff} we obtain
\begin{align}\label{poly_suff_2}
\left(p_i - \frac{1}{\delta}\right)\left(\frac{1}{p_i} - \frac{1}{\delta}\right) &= \alpha_i^2 + \frac{\epsilon_i^2}{1 + \epsilon_i} > \alpha_i^2.
\end{align}
Hence combining \eqref{poly_suff_1} and \eqref{poly_suff_2}, we obtain an equivalent condition for \eqref{poly_suff}
\begin{align}\label{poly_suff_3}
\left(1 - \frac{1}{\delta}\right)^2 > \alpha_i^2,
\end{align}
for all $\alpha_i$. As the Right Hand Side of the Eq. \eqref{poly_suff_3} is identical for all $\alpha_i$, an equivalent condition for \eqref{poly_suff} and \eqref{poly_suff_3} is
\begin{align}\label{poly_suff_4}
\left(1 - \frac{1}{\delta}\right)^2 > \max_{\lambda_i \in \{\tilde{\lambda}_2,\ldots,\tilde{\lambda}_{N^d}\}}\alpha_i^2.
\end{align}
As $\alpha_i = (\mu^2+\sigma^2)g^2\tilde{\lambda}_i^2 - 2a_0\mu g\tilde{\lambda}_i + a_0^2$ is a quadratic in $\tilde{\lambda}_i$ with a positive corfficient for the quadratic term, the maximum over an interval can be achieved only at the end points of the interval. Hence we must have
\begin{align}\label{lam_max}
\max_{\lambda_i \in \{\tilde{\lambda}_2,\ldots,\tilde{\lambda}_{N^d}\}}\alpha_i^2 = \max_{\tilde{\lambda}_2,\tilde{\lambda}_{N^d}}\alpha_i^2.
\end{align}
Hence, the equivalent condition for \eqref{poly_suff} based on Eqs. \eqref{poly_suff_4} and \eqref{lam_max} is given by
\begin{align}\label{poly_suff_eq}
\left(1 - \frac{1}{\delta}\right)^2 > \max_{\tilde{\lambda}_2,\tilde{\lambda}_{N^d}}\alpha_i^2.
\end{align}
This gives the sufficient condition for mean square exponential synchronization for the system \eqref{torus_sys_1d}. Using the condition provided in \eqref{poly_suff_eq}, we can define a margin of synchronization for a given uncertainty variance $\sigma^2$, which will quantify how vulnerable the system is to additional uncertainty leading to a desynchronized state. This margin of synchronization can be defined based on \eqref{poly_suff_eq} as,
\begin{align}\label{margin}
\rho_{SM} := 1 - \sigma^2 \left(\frac{\tilde{\lambda}_{sup}^2g^2}{\left(1 - \frac{1}{\delta}\right)^2 - \left(a_0 - \mu\tilde{\lambda}_{sup}g\right)^2}\right),
\end{align}
where $\tilde{\lambda}_{sup} = \arg\max_{\tilde{\lambda_2},\tilde{\lambda}_{N^d}}\alpha_i^2$.
\end{IEEEproof}

\section{Computationally Efficient Sufficient Condition for Synchronization}\label{section_reducedsuffcond}
In this section, we exploit  the identical nature of  network component dynamics to derive more conservative sufficient condition that the one derived in Section \ref{section_main}. The derived sufficient condition is computationally efficient and is independent of network size.
The new sufficient condition is also very insightful as it allows us to understand the tradeoff and interplay of role played by the network property, in particular the second smallest and largest eigenvalues of the nominal interconnection Laplacian, and the statistics of uncertainty in network synchronization. We start with the following definition of coefficient of dispersion.

%
\begin{definition}[Coefficient of Dispersion]
Let $\xi \in \mathbb{R}$ be a random variable with mean $\mu > 0$ and variance $\sigma^2 > 0$. Then, the coefficient of dispersion $\gamma$ is defined as
\[\gamma := \frac{\sigma^2}{\mu}\]
\end{definition}
To utilize the above definition in subsequent results we make an assumption on the system
\begin{assumption}
For all edges $(i,j)$ in the network, the mean weights assigned are positive, i.e. $\mu_{ij} > 0$ for all $(i,j)$. Furthermore, the coefficient of dispersion of each link is given by $\gamma_{ij} = \frac{\sigma_{ij}^2}{\mu_{ij}}$, and $\bar{\gamma} = \max_{\forall (i,j)}\{\gamma_{ij}\}$. This assumption simply states that the network connections are positively enforcing the coupling.
\end{assumption}
The following theorem provides a sufficiency condition for synchronization of the coupled systems based on the stability of a single modified system.

\begin{theorem}\label{main_sync}
The coupled system \eqref{coup_dyn_comp} is mean square exponentially synchronized if there exists a symmetric positive definite matrix $P > 0$ such that $\Sigma_1 - B^\top PB > 0$ and
\begin{align}\label{suff_ric}
P & =(A_0 - \lambda_{sup} GC)^\top P(A_0 - \lambda_{sup} GC) + (A_0 - \lambda_{sup} GC)^\top PB(\Sigma_1 - B^\top PB)^{-1}B^\top P(A_0 - \lambda_{sup} GC)\nonumber\\
&\quad 2\bar{\gamma}\tau\lambda_{sup} \left(C^\top G^\top PGC + C^\top G^\top PB(\Sigma_1 - B^\top PB)^{-1}B\top PGC\right)+ C^\top \Sigma_1^{-1}C + R
\end{align}
for $R > 0$, $A_0 = A - B\Sigma_1^{-1}C$ and $\lambda_{sup} \in \{\lambda_2,\lambda_N\}$, where $\lambda_N$ is the largest eigenvalue and $\lambda_2$ is the Fiedler eigenvalue, of the nominal Laplacian. Furthermore, $\tau := \frac{\lambda_{N_u}}{\lambda_{N_u} + \lambda_{2_d}}$, where $\lambda_{N_u}$ is the largest eigenvalue of the Laplacian for the purely uncertain graph $L_u$ and $\lambda_{2_d}$ is the second smallest eigenvalue of the purely deterministic Laplacian $L_d$.
\end{theorem}
\begin{IEEEproof}Please refer to the Appendix section of this paper for the proof.
\end{IEEEproof}

\subsection{Significance of $\tau$}
In Theorem \ref{main_sync}, the factor $\tau := \frac{\lambda_{N_u}}{\lambda_{N_u}+\lambda_{2_d}}$ captures the effect of location and number of uncertain links, whereas $\bar{\gamma}$ captures the effect of intensity of the randomness in the links. It is clear that $0 < \tau \leq 1$. If the number of uncertain links ($|E_U|$) is sufficiently large, the graph formed by purely deterministic edge set may become disconnected. This will imply $\lambda_{2_d} = 0$, and, $\tau = 1$. Hence, for large number of uncertain links, $\lambda_{N_u}$ is large while $\lambda_{2_d}$ is small. In contrast, if a single link is uncertain, say $E_U = \{e_{kl}\}$, then $\tau = \frac{2\mu_{kl}}{2\mu_{kl}+\lambda_{2_d}}$. Hence, for a single uncertain link, the weight of the link has a degrading effect on the synchronization margin. The location of such an uncertain link will determine the value of $\lambda_{2_d} \leq \lambda_2$, thus degrading the synchronization margin. Based upon this observation, we can rank order individual links within a graph, with respect to their degradation of the synchronization margin, on the basis of location ($\lambda_{2_d}$), mean connectivity weight ($\mu$), and the intensity of randomness given by CoD $\gamma$.

\subsection{Significance of Laplacian Eigenvalues}\label{section_Lapeigenvalues}

The condition for synchronization in Theorem \ref{main_sync} is provided in terms of both the second and the largest eigenvalues of the mean Laplacian. While the significance of the second smallest eigenvalue of the Laplacian in terms of graph connectivity is well-known in the literature, the significance of the largest eigenvalue of the Laplacian is not well documented.
The second smallest eigenvalue, $\lambda_2>0$, of the graph Laplacian indicates algebraic connectivity of the graph. We observe from Theorem \ref{main_sync}, as equation \eqref{suff_ric1} is a quadratic in $\lambda$, there exist critical values of $\lambda_2$($\lambda_N$) for the given system parameters and CoD, below(above) which synchronization is not guaranteed, respectively. Hence, critical $\lambda_2$ indicates we require a minimum degree of connectivity within the network to accomplish synchronization. To understand the significance of $\lambda_N$, we look at the complement of the graph on the same set of nodes. We know from \cite{Merris1994}, sum of largest Laplacian eigenvalue of a graph and second smallest Laplacian eigenvalue of its complemet is constant. Thus, if $\lambda_N$ is large the complementary graph has low algebraic connectivity. Thus, if we have hub nodes with high connectivity, then these nodes are sparsely connected in the complementary graph. Thus we interpret a high $\lambda_N$ indicates a high presence of densely connected hub nodes. Therefore we conclude strong robustness property in synchronization is guaranteed for close to average connectivity of nodes as compared to isolated highly connected hub nodes.

\section{Simulation Results}\label{sim_res}
In this section we provide simulations for nonlinear system synchronization based on our results presented in this paper.

\subsection{Network of Chua's Circuit Systems}
\begin{figure}[ht]
\begin{center}
\mbox{
\hspace{-0.15in}
\subfigure[]{\scalebox{0.2}{\includegraphics{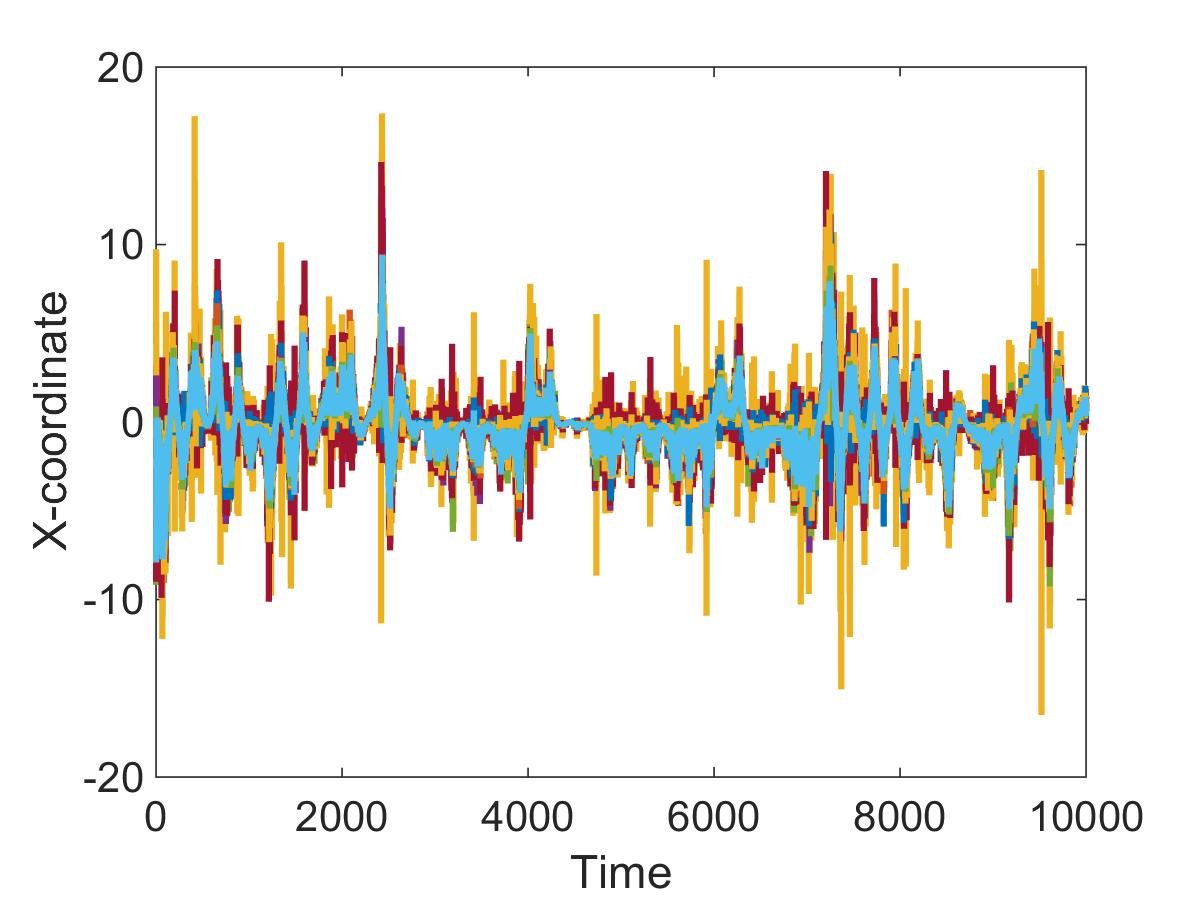}}}
\subfigure[]{\scalebox{0.2}{\includegraphics{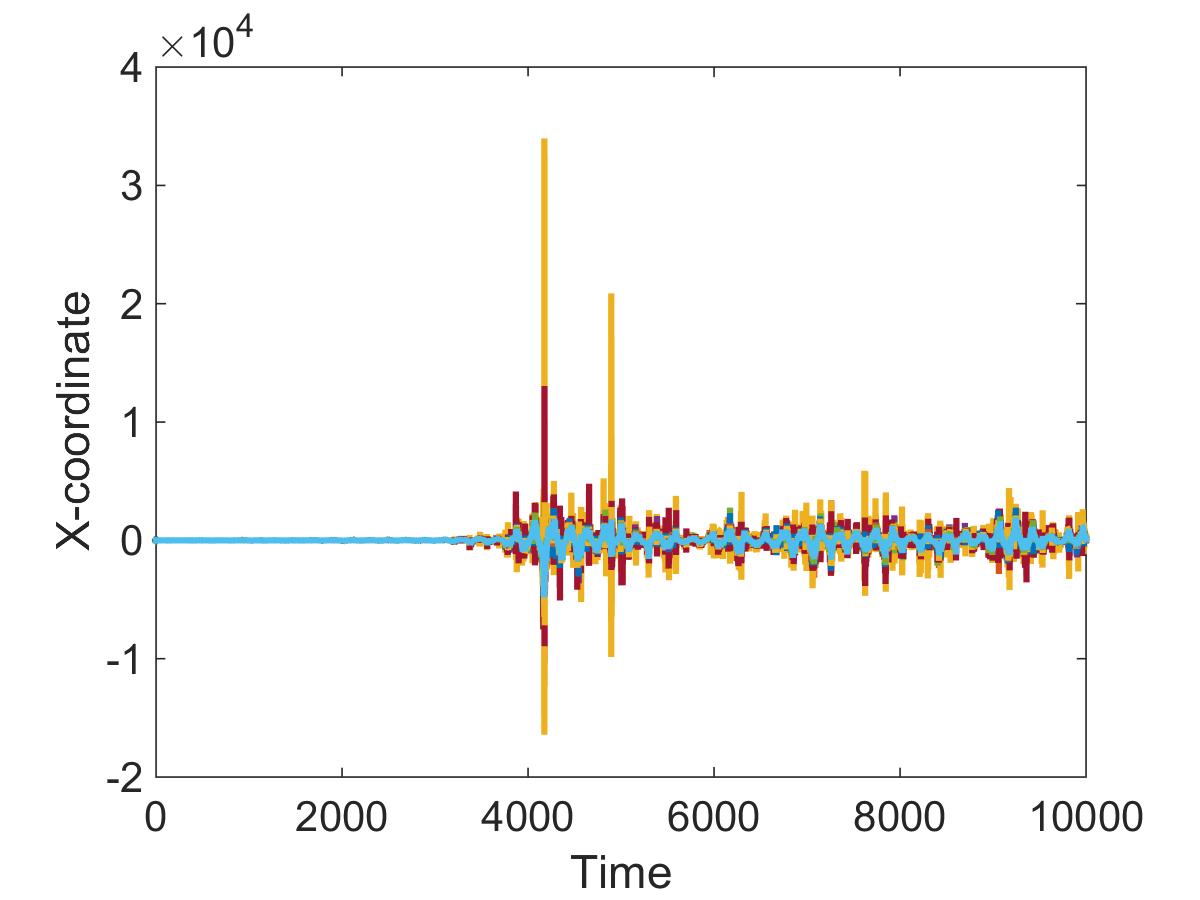}}} }
\caption{(a) X-dynamics for $\bar{\gamma}=0.9$, (b) X-dynamics for $\bar{\gamma}=1.3$}
\vspace{-0.25in}
\label{fig_chua}
\end{center}
\end{figure}

We consider network of coupled Chua's circuit systems with linear coupling and stochastic uncertainty in their interactions. The dynamics of the individual systems is given by
\[\dot x=\begin{pmatrix}0 & 7.5 & 0\\1 & -1 & 1\\ 0 & -15 & 0\end{pmatrix}x-\begin{pmatrix}7.5 \\ 0 \\ 0\end{pmatrix}\phi(y),\;\;\;y=\begin{pmatrix}1 & 0 & 0\end{pmatrix}x\]
\begin{align*}
\phi(y) = \left\{\begin{array}{cc} \epsilon y & |y| < 1 \\ (\epsilon - m_0 + m_1)y + (m_0 - m_1)\text{sgn}(y) & |y| > 1 \end{array} \right.
\end{align*}
The above system is then discretized using a zero order hold. We assume that the nonlinearity and the network interaction change only at discrete intervals and are constant during an interval. We choose the sampling time to be $T=0.01$ seconds. In \ref{fig_chua} we show the x-dynamics above and below the critical $\bar{\gamma}_c = 1.118$ below which the system should synchronize. We observe that at $\bar{\gamma} = 0.9 < \bar{\gamma}_c$ the system is synchronized. At $\bar{\gamma} = 1.3 > \bar{\gamma}_c$ the system is de synchronized.

\subsection{Spatially Periodic Systems}

We consider a simple spatially periodic system with linear dynamics $a = 1.05$, $\delta = 8$, $g=0.01$, $\mu = 1$ and $\sigma^2 = 0.01$. We take spatially periodic networks with $N=50$ agents on each dimension. We choose the number of neighbors per dimension to vary between $1$ to $25$, and the dimensions to vary between $1$ to $10$. We now plot the results for the synchronization margin as given by Theorem \ref{1Dtorus_suff_condn}.

\begin{figure}[ht]
\begin{center}
\mbox{
\hspace{-0.15in}
\subfigure[]{\scalebox{0.4}{\includegraphics{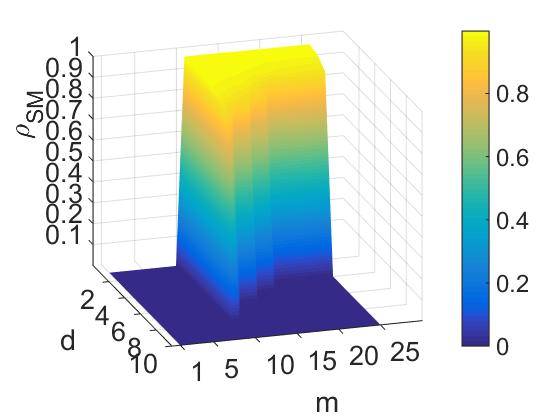}}}
\subfigure[]{\scalebox{0.18}{\includegraphics{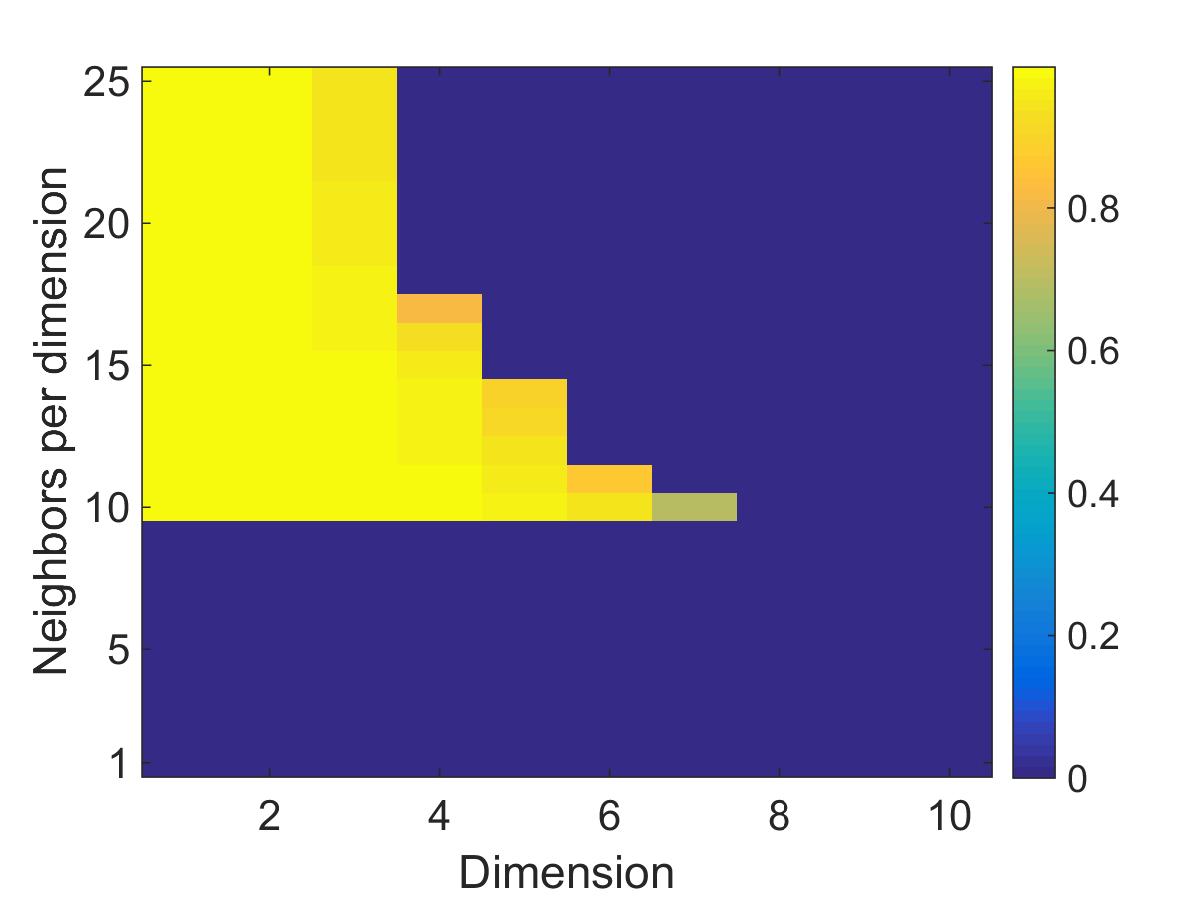}}} }
\caption{(a) 3-dimensional plot in $d-m-\rho_{SM}$ space with torus-dimension (d) on x-axis, neighbors per dimension (m) on y-axis and synchronization margin $\rho_{SM}$ on z-axis, (b) 2-dimensional plot in $d-m$ space with torus-dimension (d) on x-axis and neighbors per dimension (m). Color indicates synhcronization margin.}
\vspace{-0.25in}
\label{torus_net}
\end{center}
\end{figure}

We observe in Fig. \ref{torus_net}(a) we plot the $d-m-\rho_{SM}$ space which shows the possible synchronization margin for each value of dimension (d) and neighbors ($m$). One observes, as the number of neighbors are increased for small number of dimensions, the synchronization margin goes up. The better the connectivity, the better the margin of synchronization for smaller dimensions. As the dimension of the torus is increased, the number of optimal neighbors starts to decrease. Thus for higher dimension, it is better not to have a high number of neighbors in order to have high synchronization margin.

\section{Conclusions}\label{conc}
In this paper we study the problem of synchronization of Lur'e systems over an uncertain network. This problem is presented as a special case of the problem of stabilization of Lur'e system with parametric uncertainty. Other special case of this problem include control of Lur'e system over an uncertain network which have been previously studied by the authors. These results are used to obtain some insightful results for the problem of synchronization over uncertain networks. We conclude that the sufficient condition for mean square exponential syncronization, of the coupled dynamics, is governed by mean square exponential stability of a representative system, with multiplicative parametric uncertainty in the state matrix. This uncertainty multiplies an output feedback based on the coupling matrix, that modifies the system dynamics. The uncertainty in the reprenstative system, has a CoD twice that of the maximum CoD in the network links and its mean is a function of the eigenvalue, of the mean network Laplacian.

As the sufficient condition is based on a single representative system, it is attractive from the point of view of computational complexity for large scale networks. This sufficient condition is solved as an LMI using Schur complement, similar to deterministic Positive Real Lemma. Furthermore, these results can be used to determine the maximum amount of dispersion tolerable within the network links. As expected we conclude that, if the randomness in the network links is highly clustered then it will be more difficult to synchronize the system. Another point of interest is that the synchronization of complex nonlinear systems, depends on the largest mean Laplacian eigenvalue along with the Fiedler eigenvalue, as opposed to just one for stable or marginally stable systems achieving consensus. This indicates that while, a certain minimum connectivity needs to be present to achieve synchronization, a high density of connections among nodes might be too much information for complex nonlinear system to synchronize under uncertainty.

\section{Acknowledgements}

The authors acknowledge the support from the NSF grants ECCS 1002053, ECCS 1150405, CNS-1329915 towards this work.

\bibliographystyle{IEEEtran}
\bibliography{ref_new}


\section{Appendix}

In the appendix we provide proofs for some of the important results we prove in the paper.

\begin{IEEEproof} [{\bf Theorem \ref{main_kyp1}}]
We show the conditions in Theorem \ref{main_kyp1} are indeed sufficient by constructing an appropriate Lyapunov function that guarantees mean square exponential stability. We will prove the result in Theorem \ref{main_kyp1} for Case 1 and prove Case 2 as its dual. First, note that \eqref{unc_lure_p} holds if and only if
\begin{align} \label{unc_lure_p1}
P &= E_{\Xi(t)} \left [(A^\top (\Xi(t))PB - C^\top )(\Sigma-B^\top PB)^{-1}(C - B^\top P A(\Xi(t)))\right] \nonumber\\
& \quad + E_{\Xi(t)} \left [ A^\top (\Xi(t)) PA(\Xi(t))\right] + R_P
\end{align}
The equivalence of the two equations \eqref{unc_lure_p} and \eqref{unc_lure_p1} is observed based on \cite{lancaster} (Proposition 12.1,1). Now consider the Lyapunov function $V(x_t) = x_t^\top Px_t$. Then, the condition for the system to be mean square exponentially stable is given by
\begin{align}\label{lyap1}
E_{\Xi(t)} \left  [ V(x_{t+1}) - V(x_t) \right ] = & x^\top _t \left ( E_{\Xi} \left [A^\top (\Xi)PA(\Xi) \right] - P \right ) x_t + 2x^\top _t E_{\Xi} [A^\top (\Xi) BP]\phi(y_t, t) \nonumber\\
&+ \phi^\top (y_t,t) B^\top PB \phi(y_t,t)
\end{align}
Substituting from \eqref{unc_lure_p1} in \eqref{lyap1} and applying algebraic manipulations as adopted in \cite{haddad_bernstein_DT}, we get
\begin{align*}
 E_{\Xi(t)} \left[V(x_{t+1})\right] - V(x_t) =& - x^\top_t R_P x_t - E_{\Xi(t)}\left[\zeta^\top_t \zeta_t \right] - 2 \phi^\top (y_t,t)\left(y_t - D\phi(y_t,t)\right)
\end{align*}
where $\zeta_t = \Sigma_P^{-\frac{1}{2}}\left(B^\top PA(\Xi(t)) - C\right)x_t - \Sigma_P^{\frac{1}{2}}\phi(y_t,t)$ and $\Sigma_P =~ (\Sigma - B^\top PB)$. From condition given in Assumption \ref{nonlin_assume} we get $\phi^\top (y_t,t)\left(y_t - D\phi(y_t,t)\right)>0$, which gives us,
\[ E_{\Xi} \left  [ V(x_{t+1})) - V(x_t) \right ] < - x^\top_t R x_t < 0  \]
This implies mean square exponential stability of $x_t$ and hence Case 1 is proved. Case 2 is now the dual to Case 1 by a simple argument as shown in \cite{lure_stab_unc}.
\end{IEEEproof}

\begin{IEEEproof} [{\bf Lemma \ref{con_syn_sta}}]
From \eqref{transform1} we have
\begin{align}
\parallel \hat z_t \parallel^2 &
= \tilde x_t^\top \left(U_m\otimes I_n\right)\left(U_m^\top \otimes I_n\right)\tilde x_t = \tilde x_t^\top \left(U_mU_m^\top \otimes I_n\right)\tilde x_t \label{diff_eq1}
\end{align}
Applying $U_mU_m^\top = V_mV_m^\top - \frac{\bf{1}}{\sqrt{N}}\frac{\bf{1}^\top}{\sqrt{N}} = I_N - \frac{1}{N}\bf{1}\bf{1}^\top$ in \eqref{diff_eq1} we get
\begin{align*}
\parallel \hat z_t \parallel^2 &= \tilde x_t^\top \left(\left(I_N - \frac{1}{N}\bf{1}\bf{1}^\top\right)\otimes I_n\right)\tilde x_t= \tilde x_t^\top \left(I_{Nn} - \left(\frac{1}{\sqrt{N}}\bf{1}\otimes I_n\right)\left(\frac{1}{\sqrt{N}}\bf{1}\otimes I_n\right)^\top\right)\tilde x_t\nonumber\\
&= \tilde x_t^\top\tilde x_t - \bar x_t^\top\bar x_t = \frac{1}{2N}\sum_{i=1}^N\sum_{j\neq i,j=1}^N \left(x_t^i - x_t^j\right)^\top\left(x_t^i - x_t^j\right)
\end{align*}
Now, mean square exponential stability of \eqref {coup_dyn_comp2} implies there exists $K > 0$ and $0 < \beta < 1$ such that
\begin{align*}
E_{\Xi} \parallel \hat z_{t} \parallel^2 &\le   K  {\beta}^t \parallel \hat z_{0}  \parallel^2,  \\
E_{\Xi} \sum_{k=1}^N\sum_{j\neq k, j=1}^N \parallel x^k_t - x^j_t \parallel^2 &\le  K  {\beta}^t  \sum_{k=1}^N\sum_{j\neq k,j=1}^N \parallel x^k_0 - x^j_0  \parallel^2,  \\
\Rightarrow \sum_{k=1}^N\sum_{j\neq k, j=1}^N E_{\Xi}\parallel x^k_t -x^j_t \parallel^2 &\le  K  {\beta}^t  \sum_{k=1}^N\sum_{j\neq k,j=1}^N \parallel x^k_0 - x^j_0  \parallel^2,
\end{align*}
This gives us the result,
\begin{align*}
E_{\Xi}  \parallel x^k_t - x^l_t \parallel^2 \le \bar{K}(\tilde{e}_0)\beta^t \parallel x^k_0 - x^l_0 \parallel^2.
\end{align*}
where $\bar K(\tilde{e}_0) := K \left ( 1 + \frac{\sum_{i=1, i\neq k}^N\sum_{j=1, j\neq i}^N \parallel x^i_0 - x^j_0 \parallel^2}{\parallel x^k_0 - x^l_0 \parallel^2}  \right )$.
\end{IEEEproof}
\begin{IEEEproof} [{\bf Theorem \ref{main_sync}}]
We know mean square exponential synchronization is guaranteed by conditions in Lemma \ref{suf_syn1}. Consider $\mathcal{P} = I_{N-1}\otimes P$ where $P > 0$ is a symmetric positive definite matrix that satisfies $\Sigma_1 - B^\top PB > 0$. This gives us $\hat{\Sigma_1} - \hat{B}^\top \mathcal{P}\hat{B} > 0$. Using this we write \eqref{ric_syn} as
\begin{align} \label{ric1}
& I_{N-1}\otimes P >  (\hat{A}_0 - \Lambda_m\otimes GC)^\top(I_{N-1}\otimes P)(\hat{A}- \Lambda_m\otimes GC) + \sum_{\mathcal{I}} \sigma_{\alpha_k}^2A_{\alpha_k}^\top (I_{N-1}\otimes P) A_{\alpha_k}\nonumber\\
& +(\hat{A}_0-\Lambda_m\otimes GC)^\top (I_{N-1}\otimes P)\hat{B}\left(\hat \Sigma_1 - \hat B^\top(I_{N-1}\otimes P)\hat B\right)^{-1}\hat B^\top(I_{N-1}\otimes P)(\hat{A}_0-\Lambda_m\otimes GC) \nonumber \\
& +\sum_{\mathcal{I}} {\sigma}_{\alpha_k}^2 A_{\alpha_k}^\top (I_{N-1}\otimes P)\hat B \left(\hat \Sigma_1-\hat B^\top (I_{N-1}\otimes P)\hat B\right)^{-1}\hat B^\top (I_{N-1}\otimes P) A_{\alpha_k}\nonumber\\
&+ I_{N-1}\otimes C^\top \Sigma_1^{-1}C
\end{align}
Since $A_{\alpha_k} = \hat{\ell}_{ij}\hat{\ell}_{ij}^\top \otimes GC$ we can write \eqref{ric1} as
\begin{align} \label{ric2}
I_{N-1}\otimes P &>  \left[A_0 - \lambda_j GC\right]^\top(I_{N-1}\otimes P)\left[A_0 - \lambda_j GC\right]+ I_{N-1}\otimes C^\top\Sigma_1^{-1}C\nonumber\\
& +\left[A_0 - \lambda_j GC\right]^\top\left(I_{N-1}\otimes (PB\left(\Sigma_1 - B^\top PB\right)^{-1}B^\top P)\right)\left[A_0 - \lambda_j GC\right] \nonumber \\
& +\sum_{\mathcal{I}} {\sigma}_{\alpha_k}^2 (\hat{\ell}_{\alpha_k}\hat{\ell}_{\alpha_k}^\top \otimes GC)^\top \left(I_{N-1}\otimes(PB\left(\Sigma_1 - B^\top PB\right)^{-1}B^\top P)\right)(\hat{\ell}_{\alpha_k}\hat{\ell}_{\alpha_k}^\top \otimes GC) \nonumber \\
& + \sum_{\mathcal{I}} \sigma_{\alpha_k}^2(\hat{\ell}_{\alpha_k}\hat{\ell}_{\alpha_k}^\top  \otimes GC)^\top (I_{N-1}\otimes P) (\hat{\ell}_{\alpha_k}\hat{\ell}_{\alpha_k}^\top \otimes GC)
\end{align}
where $\left[A_0 - \lambda_j GC\right] = (\hat{A}_0-\Lambda_m\otimes GC)$. Inequality \eqref{ric2} can be further simplified as
\begin{align} \label{ric3}
I_{N-1}\otimes P &>  \left[A_0 - \lambda_j GC\right]^\top(I_{N-1}\otimes P)\left[A_0 - \lambda_j GC\right]+ 2\sum_{\mathcal{I}} \sigma_{\alpha_k}^2\hat{\ell}_{\alpha_k}\hat{\ell}_{\alpha_k}^\top \otimes C^\top G^\top PGC\nonumber\\
& +\left[A_0 - \lambda_j GC\right]^\top \left(I_{N-1}\otimes (PB\left(\Sigma_1 - B^\top PB\right)^{-1}B^\top P)\right)\left[A_0 - \lambda_j GC\right] \nonumber \\
& +2\sum_{\mathcal{I}} {\sigma}_{\alpha_k}^2 \hat{\ell}_{\alpha_k}\hat{\ell}_{\alpha_k}^\top \otimes \left(C^\top G^\top PB\left(\Sigma_1 - B^\top PB\right)^{-1}B^\top PGC\right)+ I_{N-1}\otimes C^\top \Sigma_1^{-1}C
\end{align}
We know that
\begin{align}\label{co_dis1}
\sum_{\mathcal{I}}\sigma_{\alpha_k}^2\hat{\ell}_{\alpha_k}\hat{\ell}_{\alpha_k}^\top
= \sum_{\mathcal{I}}\gamma_{\alpha_k}\mu_{\alpha_k}\hat{\ell}_{\alpha_k}\hat{\ell}_{\alpha_k}^\top
\leq \bar{\gamma}\sum_{\mathcal{I}}\mu_{\alpha_k}\hat{\ell}_{\alpha_k}\hat{\ell}_{\alpha_k}^\top = \bar{\gamma}U_m^\top L_uU_m
\end{align}
We know that $L_m = L_u + L_d$. Thus if there exists $\tau \leq 1$ such that $L_u \leq \tau L_m$ then we must have $\frac{1-\tau}{\tau}L_u \leq L_d$. This is true if
\[ \left(\frac{1-\tau}{\tau}\right)\lambda_{N_u} \leq \lambda_{2_d}, \]
where $\lambda_{N_u}$, is the largest eigenvalue of the Laplacian $L_u$ and $\lambda_{2_d}$, is the second smallest eigenvalue of the Laplacian $L_d$. We now choose $\tau = \frac{\lambda_{N_u}}{\lambda_{N_u}+\lambda_{2_d}}$ and applying $L_u \leq \tau L_m$ to \eqref{co_dis1} we obtain,
\begin{align}\label{co_dis}
\sum_{\mathcal{I}}\sigma_{\alpha_k}^2\hat{\ell}_{\alpha_k}\hat{\ell}_{\alpha_k}^\top  \leq \bar{\gamma}U_m^\top (\tau L_m)U_m = \bar{\gamma}\tau\hat{\Lambda}_m
\end{align}
Now, substituting \eqref{co_dis} in \eqref{ric3} a sufficient condition for inequality \eqref{ric3} to hold is given by
\begin{align} \label{ric4}
I_{N-1}\otimes P &>  \left[A_0 - \lambda_j GC\right]^\top (I_{N-1}\otimes P)\left[A_0 - \lambda_j GC\right]+ I_{N-1}\otimes C^\top \Sigma_1^{-1}C\nonumber\\
& \quad+\left[A_0 - \lambda_j GC\right]^\top \left(I_{N-1}\otimes (PB\left(\Sigma_1 - B^\top PB\right)^{-1}B^\top P)\right)\left[A_0 - \lambda_j GC\right] \nonumber \\
& \quad+2\bar{\gamma}\tau\hat{\Lambda}_m \otimes \left(C^\top G^\top \left(P+PB\left(\Sigma_1 - B^\top PB\right)^{-1}B^\top P\right)GC\right)
\end{align}
Equation \eqref{ric4} is essentially a block diagonal equation which gives the sufficient condition for mean square exponential synchronization to be
\begin{align} \label{ric5}
P &>  (A_0 - \lambda_j GC)^\top P(A_0 - \lambda_jGC) +(A_0 - \lambda_j GC)^\top PB\left(\Sigma_1 - B^\top PB\right)^{-1}B^\top P(A_0 - \lambda_j GC) \nonumber \\
& \quad + 2\bar{\gamma}\tau\lambda_j C^\top G^\top PGC + 2\bar{\gamma}\tau\lambda_jC^\top G^\top PB\left(\Sigma_1 - B^\top PB\right)^{-1}B^\top PGC + C^\top \Sigma_1^{-1}C
\end{align}
for all non-zero eigenvalues $\lambda_j$ of $\hat{\Lambda}_m$. Since \eqref{ric5} is a quadratic in the eigenavlues $\lambda_j$, it is sufficienct to study is the equations holds true for the extreme values of the set given by $\lambda_2$ and $\lambda_N$. This is easily seen by the following argument. Using Schur complement we can equivalently write \eqref{ric5} for a given $\lambda_j$ and $C_1 = \sqrt{2\bar{\gamma}\tau}C$ as an LMI given by
\begin{align} \label{lmi_compact}
M_1 + \lambda_j M_2 > 0
\end{align}
where
\begin{align*}
M_1 = \left[\begin{array}{ccccc}
P - C^\top \Sigma_1^{-1}C & A_0^\top P & A_0^\top PB & 0 & 0\\
PA_0 & P & 0 & 0 & 0 \\
B^\top PA_0 & 0 & \Sigma_1 - B^\top PB & 0 & 0 \\
0 & 0 & 0 & 0 & 0 \\
0 & 0 & 0 & 0 & 0
\end{array}\right],
\end{align*}
\begin{align*}
\setlength{\arraycolsep}{2pt}
M_2 = \begin{bmatrix}
0 & -C^\top G^\top P & -C^\top G^\top PB & C_1^\top G^\top P & C_1^\top G^\top PB\\
-PGC & 0 & 0 & 0 & 0 \\
-B^\top PGC & 0 & 0 & 0 & 0 \\
PGC_1 & 0 & 0 & P & 0 \\
B^\top PGC_1 & 0 & 0 & 0 & \Sigma_1 - B^\top PB
\end{bmatrix}.
\end{align*}
Since this is a convex constraint in $\lambda$, if it is satisfied for any values of $\lambda_i,\lambda_j \in \{\lambda_2,\ldots,\lambda_N\}$, then \eqref{lmi_compact} is true for any $\lambda = s\lambda_i + (1-s)\lambda_j$ for all $s \in [0,1]$. Thus if we require \eqref{ric5} to hold for all eigenvalues of the mean Laplacian matrix, then it must hold for the extreme points of the set, i.e. $\lambda_{sup} \in \{\lambda_2,\lambda_N\}$. This proves the result.
\end{IEEEproof}

\end{document}